\documentclass[12pt]{amsart}
\makeatletter

\newcommand{\kk}{\kappa}

\newcommand{\ff}{\varphi}
\newcommand{\CC}{\mathbb C}

\newcommand{\ZZ}{\mathbb Z}

\newcommand{\LyX}{L\kern-.1667em\lower.25em\hbox{Y}\kern-.125emX\spacefactor1000}

\theoremstyle{plain}    
\newtheorem{thm}{Theorem}[section]
\numberwithin{equation}{section} 
\numberwithin{figure}{section} 
\theoremstyle{plain}    
\newtheorem*{thm*}{Theorem} 
\theoremstyle{plain}    
\newtheorem{cor}[thm]{Corollary} 
\theoremstyle{plain}    
\newtheorem{prop}[thm]{Proposition} 
\theoremstyle{remark}
\newtheorem{rems}[thm]{Remarks}

\theoremstyle{definition}
\newtheorem{defn}[thm]{Definition}
\makeatother

\begin{document}

\title[Operator-valued distributions I]{Operator-valued
distributions.\\ 
{\small I. Characterizations of Freeness}}
\
\author[A. Nica]{Alexandru Nica{*} }

\address{Department of Pure Mathematics, University of Waterloo, Waterloo, Ontario,
N2L 3G1, Canada}

\email{anica@math.uwaterloo.ca}

\thanks{{*} Research supported by a grant of NSERC, Canada.}

\author[D. Shlyakhtenko]{Dimitri Shlyakhtenko\dag}

\address{Department of Mathematics, UCLA, Los Angeles, CA 90095, USA}

\email{shlyakht@math.ucla.edu}

\thanks{\dag Research supported by an NSF postdoctoral fellowship, by NSF grant
DMS-0102332, and by a Sloan Foundation Fellowship}

\author[R. Speicher]{Roland Speicher{*}}

\address{Department of Mathematics and Statistics, Queen's University, Kingston, Ontario
K7L 3N6, Canada}

\email{speicher@mast.queensu.ca}

\begin{abstract}
Let $M$ be a $B$-probability space. Assume that $B$ itself is a $D$-probability
space; then $M$ can be viewed as $D$-probability space as well. Let $X\in M$. We
look at the question of relating the properties of $X$ as $B$-valued
random variable to its properties as $D$-valued random variable. We characterize
freeness of $X$ from $B$ with amalgamation over $D$: (a) in terms of a 
certain factorization
condition linking the $B$-valued and $D$-valued cumulants of $X$, and 
(b) for $D$ finite-dimensional, in terms of linking
the $B$-valued and the $D$-valued Fisher information of $X$. We give an application
to random matrices. For the second characterization we
derive a new operator-valued description of the conjugate variable and introduce an
operator-valued version of the liberation gradient.
\end{abstract}
\maketitle

\section{Introduction}

Free probability theory is a non-commutative probability theory where
the classical concept of independence is replaced by the notion of
"freeness". This theory, due to Voiculescu, was introduced as a tool
for investigating the structure of von Neumann algebras arising from
free product constructions. This programme has been very succesful
and has yielded a wealth of new and unexpected results about this class
of von Neumann algebras.

 From the very beginning, Voiculescu introduced also an operator-valued
version of freeness -- where the role of the "constants" $\CC$ is taken
over by an arbitrary fixed algebra $B$ and where the states are
replaced by conditional expectations onto $B$. 
This more general frame enlarges the domain of applicability of free
probability techniques in a tremendous way. 
Of course, this wider domain of applicability is compensated by the fact
that it is harder to obtain results in the operator-valued case.
However, quite astonishingly, a lot of the
scalar-valued theory, in particular its combinatorial description resting
on the notion of cumulants, can be transfered to the operator-valued context.

A systematic exploration about how much of the
scalar-valued results can be generalized to the operator-valued setting is
still lacking. However, instead of pursuing 
such a generalization for its own sake, we 
will develop here those aspects of such a programme which are related to 
one of the most exciting
possibilities available in the operator-valued framework: 
the possibility of
switching between two different algebras of "scalars". 
Namely if the algebra $A$ is simultaneously a $B$-probability space and
a $D$-probability space for some subalgebras $D$ and $B$ of $A$, then an
element $X$ in $A$ is at the same time a $B$- and a $D$-valued random
variable, and one can ask how these different points of view towards $X$
are related.
Let us assume that $D\subset B$. Then, in principle,
the $B$-valued distribution of $X$ determines also the $D$-valued distribution
of $X$. However, this connection is 
similar in complexity to saying that the
entries of a matrix determine its eigenvalues.
What we are looking for are treatable and interesting special
cases where something more explicit can be said.
 
Note that such questions are also of practical relevance, since in concrete
cases one might be interested in the $D$-valued distribution of $X$, however,
direct arguments only give information about its $B$-valued distribution. Then
it is of great importance to have general theorems about closing the gap between
$B$ and $D$.

In this article we will consider the most fundamental special case, namely 
when $X$ is free from $B$ over $D$. We will show that
this freeness condition can equivalently be characterized in terms of 
cumulants and in terms
of free Fisher information.

In Section 2, we will recall some preliminaries about operator-valued
free probability theory, in particular the concepts of cumulants and
canonical random variables.

Section 3 deals with our first main result: freeness from a 
subalgebra $B$ can be characterized by a factorization property of the 
$B$-valued cumulants. We also give an interesting application of this
circle of ideas to Gaussian random band matrices.

The rest of the paper deals with our second characterization of freeness
from $B$ over $D$, in the case of finite-dimensional $D$. The question
which we address is whether equality of the free Fisher information with
respect to $D$ and the free Fisher information with respect to $B$ 
implies freeness from $B$ over $D$. A more suggestive form of this question
might be the following version:
If we know that the free entropy of a random variable $X$ conditioned on
random variables $Y$ and $Z$ is the same as the free entropy of $X$ conditioned
on $Z$, does this imply that $X$ and $Y$ are conditionally free over $Z$?
(This is related to the previous questions by setting $B=W^*(Y,Z)$ and
$D=W^*(Z)$.)

For the case $D=\CC$ the above question was solved in the affirmative by 
Voiculescu \cite{V6},
but the general case remained open. We are able to extend the affirmative solution
of Voiculescu to the case where $D$ is finite-dimensional. 

In Section 4, we recall the concepts of
relative Fisher information and conjugate variables. Our main result in
this section is a reformulation
of the characterizing equations for conjugate variables in terms of 
operator-valued cumulants. 

Our main tool for addressing the above mentioned problem on Fisher informations
is an operator-valued generalization of the liberation gradient, which we will
define in Section 5. Again, we give an interesting reformulation of it in
terms of cumulants and, for the case of finite-dimensional $D$, we prove a 
relation between the liberation gradient and corresponding
conjugate variables. 

These results are used in Section 6 for proving, in the case of finite-dimensional
$D$, our
second main characterization for freeness from a subalgebra.
We will also reformulate this result as a minimization result for free Fisher
information or a maximization result for free entropy. Finally, we treat as
a concrete example of such a maximization result the case of $R$-cyclic elements.

A preliminary version of sections 2 and 3 of
this paper has appeared as MSRI-preprint 2001-001.
The authors gratefully acknowledge the hospitality of MSRI during its
2000-2001 program in operator algebras.

\section{Preliminaries}

\subsection{$B$-probability space}
Let \( B \) be a unital algebra. Recall that a \emph{\( B \)-probability space}
$(M,\,E:M\to B)$ (see e.g. \cite{DVV:book}, \cite{speicher:thesis}) is
a pair consisting of an algebra \( M \) containing \( B \) as a unital subalgebra,
and a \emph{conditional expectation}~\( E:M\to B \). In other words, \( E:M\to B \)
is unital and \( B \)-bilinear:
\[
E(1)=1,\qquad E(bmb')=bE(m)b',\quad \forall b,b'\in B,\forall m\in M.\]
Elements of \( M \) are called \( B \)-valued random variables. 

\subsection{Multiplicative functions.}
Recall that a \emph{$B$-balanced map} \( \langle \cdots \rangle :\bigcup M^{n}\to B \)
is a \( \mathbb {C} \)-multilinear map, (i.e., a sequence of maps \( M^{n}\ni m_{1},\dots ,m_{n}\mapsto \langle m_{1},\dots ,m_{n}\rangle \in B \))
satisfying the \( B \)-linearity conditions
\begin{eqnarray*}
\langle bm_{1},\dots ,m_{n}\rangle  & = & b\langle m_{1},\dots ,m_{n}\rangle \\
\langle m_{1},\dots ,m_{n}b\rangle  & = & \langle m_{1},\dots ,m_{n}\rangle b\\
\langle m_{1},\dots ,m_{k}b,m_{k+1},\dots ,m_{n}\rangle  & = & \langle m_{1},\dots ,m_{k},bm_{k+1},\dots ,m_{n}\rangle .
\end{eqnarray*}
(Here \( m_{i}\in M \), \( b\in B \)). 

Given a non-crossing partition \( \pi \in NC(n) \) and an arbitrary $B$-balanced
map \( \langle \cdots \rangle  \), we can construct a corresponding \emph{multiplicative
map} or 
\emph{bracketing}, denoted by
$\langle\cdots\rangle_\pi $, which is a
 map $M^{n}\to B$ and is  
defined recursively by
\begin{align*}
\langle m_{1},\dots ,m_{k}\rangle_{1_k}  & =  \langle m_{1},\dots ,m_{k}\rangle \\
\langle m_{1},\dots ,m_{k}\rangle_{\pi \sqcup \rho}   & =  
\langle m_{1},\dots ,m_{p}\rangle_\pi \cdot \langle 
m_{p+1},\dots ,m_{k}\rangle_\rho \\
\langle m_{1},\dots ,m_{k}\rangle_{\textrm{ins}(p,\rho \to \pi )} & =  
\langle m_{1},\dots ,m_{p}\langle m_{p+1},\dots ,m_{p+q}\rangle_\rho ,
m_{p+q+1},\dots ,m_{k}\rangle_\pi .
\end{align*}
Here \( 1_{k} \) denotes the partition with the sole class \( \{1,\dots ,k\} \),
\( \pi \sqcup \rho  \) denotes disjoint union (with the equivalence classes
of \( \rho  \) placed after those of \( \pi  \)), and \( \textrm{ins}(p,\rho \to \pi ) \)
denotes the partition obtained from \( \pi  \) by inserting the partition \( \rho  \)
after the \( p \)-th element of the set on which \( \pi  \) determines a partition. 

In other words, each partition \( \pi  \) is interpreted as a recipe for placing
brackets \( \langle \cdots \rangle  \), and \(\langle \cdots \rangle_\pi  \)
is the value of the resulting expression.

\subsection{Moments and \protect\( R\protect \)-transform.}

The \( B \)-probability space structure of the algebra \( M \) gives rise
to one example of such a multiplicative map, namely, the \emph{moments map}
\[
\mu^B:\bigcup M^{n}\to B,\]
given by
\begin{equation}
\mu^B(m_{1},\dots ,m_{n})=E(m_{1}\cdots m_{n}).
\end{equation}
The reason for the name is that, having fixed \( B \)-random variables \( X_{1},\dots ,X_{p}\in M \),
the following values of $\mu^B$,
$$\mu^B(b_{0}X_{i_{1}}b_{1},\dots ,X_{i_{n-1}}b_{n-1},X_{i_{n}}b_{n})
=E(b_{0}X_{i_{1}}b_{1}\cdots X_{i_{n-1}}b_{n-1}X_{i_{n}}b_{n})
$$
are called \emph{\( B \)-valued moments} of the family \( X_{1},\dots ,X_{n} \). 

In \cite{speicher:thesis} and \cite{dvv:amalg} the notion of \( B \)-valued
\( R \)-transform was introduced (we follow the combinatorial approach of \cite{speicher:thesis}, see also \cite{cumulants}).
Like the map $\mu^B$, the \( R \)-transform map is a multiplicative
map
$$
\kk^B:\bigcup M^{n}\to B.$$
The following combinatorial formula ("moment-cumulant formula")
actually determines $\kk^B$
uniquely:
\begin{align}
\label{moment-cumulant}
\mu^B(m_1,\dots,m_n)&=
E(m_{1}\cdots m_{n})\notag \\
&=\sum \textrm{all possible bracketings involving }
\kk^B(\cdots)  \\
&=\sum _{\pi \in NC(n)}\kk_\pi^B(m_{1},\dots ,m_{n}). \notag 
\end{align}
The uniqueness of the definition can be easily seen by observing that the right-hand
side of the equation above involves $\kk^B(m_{1},\dots ,m_{n})$ and that
the rest of the terms are products of factors of smaller order (i.e., restrictions
of $\kk^B$ to \( M^{k} \), \( k<n \)). 

It is important to note that $\mu^B$ determines $\kk^B$
and vice-versa. Moreover, the value of $\mu^B|_{M^{n}} \) depends only
on $\kk^B|_{M\cup \cdots \cup M^{n}} \), and vice-versa.

\subsection{Moment and cumulant series.}

Suppose that we are given 
a $B$-balanced function \( \langle \cdots \rangle :\bigcup M^{k}\to B \),
and that on the other hand we only want to focus our attention to a given
family
\( X_{1},\dots ,X_{n}\in M \) of elements of $M$.
Then what we have to look at is the family of multilinear maps
\[
M^{\langle \cdots \rangle }_{i_{1},\dots ,i_{k}}:B^{k-1}\to B,\quad i_{1},\dots ,i_{k}\in \{1,\dots ,n\}\]
defined by
\[
M^{\langle \cdots \rangle }_{i_{1},\dots ,i_{k}}(b_{1},\dots ,b_{k-1})=\langle X_{i_{1}}b_{1},\cdots ,b_{k-1}X_{i_{k}}\rangle .\]
In the particular examples above, we get the \emph{moment series} 
of \( X_{1},\dots ,X_{n} \),
\begin{equation}
\mu ^{X_{1},\dots ,X_{n}}_{i_{1},\dots ,i_{k}}:=M^{\mu}_{i_{1},\dots ,i_{k}}
\end{equation}
and the \emph{cumulant series},
\begin{equation}
k_{i_{1},\dots ,i_{k}}^{X_{1},\dots ,X_{n}}:=M_{i_{1},\dots ,i_{k}}^{\kk}.
\end{equation}
We will sometimes write
\[
k_{B;i_{1},\dots ,i_{k}}^{X_{1},\dots ,X_{n}}\]
to emphasize that the series is valued in \( B \).

\subsection{Freeness with amalgamation.}

Let \( M_{1},M_{2}\subset M \) be two subalgebras, each containing \( B \).
Freeness of $M_1,M_2$ with amalgamation over $B$ is defined in the same
way as the usual (scalar-valued) freeness, one just has to replace the
state $\tau$ by the conditional expectation $E:M\to B$. We refer to \cite{DVV:book}
for details. 
The importance of the $B$-valued \( R \)-transform in the context of
freeness with amalgamation over $B$  
is apparent from the following theorem
(\cite{speicher:thesis}, see also \cite{cumulants-nica}):

\begin{thm*}
\label{mixed}
Let \( S_{1},S_{2} \) be two subsets of \( M \). Let \( M_{i} \) be the algebra
generated by \( S_{i} \) and \( B \), \( i=1,2. \) Then \( M_{1} \) and
\( M_{2} \) are free with amalgamation over \( B \) iff whenever \( X_{1},\dots ,X_{n}\in S_{1}\cup S_{2} \),
\[
\kk^B(X_{1},\dots ,X_{n})=0\]
unless either all \( X_{1},\dots ,X_{n}\in S_{1} \), or all \( X_{1},\dots ,X_{n}\in S_{2} \).
\end{thm*}

Note that the above theorem makes a statement only if $n\geq 2$.

\subsection{Canonical random variables.\label{sec:canonical}}

Let \( X_{1},\dots ,X_{n}\in M \) be fixed. 

Then by \cite{dvv:amalg} there exists a \( B \)-probability space \( (\mathcal{F},E_{B}:\mathcal{F}\to B) \),
elements \( \lambda _{1}^{*},\dots ,\lambda _{n}^{*}\in \mathcal{F} \), and
elements \( \lambda _{p}^{k} \) ($1\leq p\leq n, k\geq 0$)
satisfying the following properties:\renewcommand{\labelenumi}{(\roman{enumi})}\renewcommand{\theenumi}{\arabic{section}.\arabic{subsection}(\roman{enumi})}

\begin{enumerate}
\item \label{relation:lambdas}\( \lambda _{j_{1}}^{*}b_{1}\lambda _{j_{2}}^{*}b_{2}\dots \lambda _{j_{k}}^{*}b_{k}\lambda _{j}^{k}=k_{j_{1},\dots ,j_{k},j}^{X_{1},\dots ,X_{n}}(b_{1},\dots ,b_{k}) \),
\( b_{1},\dots ,b_{k}\in B \); 
\item \label{relation:E}Let \( w=b_{0}a_{1}b_{1}a_{2}b_{2}\dots a_{n}b_{n} \), where
\( b_{i}\in B \), \( a_{i}=\lambda ^{k}_{p} \) or \( a_{i}=\lambda _{p}^{*} \).
Then \( E_{B}(w)=0 \) unless \( w \) can be reduced to an element of \( B \)
using relation \ref{relation:lambdas}.
\end{enumerate}
The construction of the elements $\lambda_p^k$ puts in particular
 \( \lambda _{p}^{0}=E_{B}(X_{p})\in B \). We should
clarify that for each \( k \), \( \lambda _{p}^{k} \) is just a formal variable,
and we do not assume any relations between \( \{\lambda _{k}^{p}\}_{k,p} \);
for example, \( \lambda _{p}^{k} \) is not the \( k \)-th power of \( \lambda _{p}^{1} \).

It is not hard to show that the properties listed above determine the restriction
of \( E_{B} \) to the algebra generated by 
\( \{\lambda_p ^{*}\}_p \) and \( \{\lambda _{p}^{k}\}_{p,k} \), and $B$. 

Let
\[
Y_{j}=\lambda _{j}^{*}+\sum _{k\geq 0}\lambda _{j}^{k}.\]
(This series is formal; however, \( Y_{1},\dots ,Y_{n} \) have moments, since
each such moment involves only a finite number of terms from the series defining
\( Y_{j} \)).

The point of this construction is that we have:
\[
k_{i_{1},\dots ,i_{k}}^{X_{1},\dots ,X_{n}}=
k_{i_{1},\dots ,i_{k}}^{Y_{1},\dots ,Y_{n}}\quad \textrm{and }\quad
\mu _{i_{1},\dots ,i_{k}}^{X_{1},\dots ,X_{n}}=
\mu _{i_{1},\dots ,i_{k}}^{Y_{1},\dots ,Y_{n}}.\]
In other words, given a cumulant series, \( Y_{1},\dots ,Y_{n} \) is an explicit
family of \( B \)-valued random variables, whose cumulant series is equal to
the one given.

\section{Freeness from a subalgebra and factorization of cumulants.}

\subsection{\protect\( D\protect \)-cumulants vs. \protect\( B\protect \)-cumulants.}

Let now \( D\subset B \) be a unital subalgebra, and let \( F:B\to D \) be
a conditional expectation. If \( (M,E:M\to B) \) is a \( B \)-probability
space, then \( (M,F\circ E:M\to D) \) is a \( D \)-probability space. 

\begin{thm}
\label{restriction}
Let \( X_{1},\dots ,X_{n}\in M \). Assume that the \( B \)-valued cumulants
of \( X_{1},\dots ,X_{n} \) satisfy
\begin{equation}
k^{X_{1},\dots ,X_{n}}_{B;i_{1},\dots ,i_{k}}(d_{1},
\dots ,d_{k-1})\in D,\quad \forall k\geq 1,\,\forall d_{1},\dots ,d_{k-1}\in D.
\end{equation}
Then the \( D \)-valued cumulants of \( X_{1},\dots ,X_{n} \) are given by
the restrictions of the \( B \)-valued cumulants:
\begin{equation}
k^{X_{1},\dots ,X_{n}}_{D;i_{1},\dots ,i_{k}}(d_{1},
\dots ,d_{k-1})=k^{X_{1},\dots ,X_{n}}_{B;i_{1},\dots ,
i_{k}}(d_{1},\dots ,d_{k-1}),\quad \forall k\geq 1,\,\forall d_{1},\dots ,d_{k-1}\in D.
\end{equation}

\end{thm}
\begin{proof}
Let \( N \) be the algebra generated by \( D \) and \( X_{1},\dots ,X_{n} \).
The condition on cumulants implies that
\[
\kk^B|_{\bigcup N^{p}}\]
is valued in \( D \). It follows from the moment-cumulant 
formula (\ref{moment-cumulant})
that $\mu^B|_{N^{p}}$
is valued in \( D \), and hence (by the simple observation that in general
$\mu^D=F(\mu^B)$) that
\[
\mu^D|_{N^{p}}=
\mu^B|_{N^{p}}=
\sum _{\pi \in NC(p)}\kk^B_\pi|_{\bigcup N^{p}}.\]
Since the moment-cumulant formula determines $\kk^D|_{\bigcup N^{p}}$,
it follows that
\[
\kk^D|_{\bigcup N^{p}}=\kk^B|_{\bigcup N^{p}}.\]

\end{proof}
We record an equivalent formulation of the theorem above (which was implicit
in the proof):

\begin{thm}
\label{thrm:restricRtransformDvalued}Let \( N\subset M \) be a subalgebra,
containing \( D \). Assume that $\kk^B|_{\bigcup N^{p}}$ is
valued in \( D \). Then
\begin{equation}
\kk^D|_{\bigcup N^{p}}=\kk^B|_{\bigcup N^{p}}.
\end{equation}

\end{thm}
In general, in the absence of the condition that \( k^{X_{1},\dots ,X_{n}}_{B} \)
restricted to \( \cup D^{p} \) is valued in \( D \), the expression of \( D \)-valued
cumulants of \( X_{1},\dots ,X_{n} \) in terms of the \( B \)-valued cumulants
is quite complicated. Note, for example, that if \( X \) is a \( B \)-valued
random variable, and \( b\in B \), then the \( B \)-valued cumulant series
of \( bX \) are very easy to describe. On the other hand, the \( D \)-valued
cumulant series of \( bX \) can have a very complicated expression in terms
of the \( D \)-valued cumulant series of \( X \) and \( b \).

In spite of its apparent simplicity, Theorem \ref{restriction} has non-trivial
applications. One kind of application appears, e.g., in the following
type of situation: We are interested in the $D$-valued cumulants for
a certain $D$, but there is no nice general formula for calculating
$D$-valued cumulants. However, we can find a larger algebra $B$, containing
$D$, where there is a nice formula for $B$-valued cumulants. Then
Theorem \ref{restriction} serves us with special situations when the
desired $D$-valued cumulants can nevertheless be computed. Situations
like this occur, for example, in the context of R-cyclic elements, as
considered in \cite{NSS2}.

The sufficient condition in the theorem above is actually quite close to being
necessary in the case that the conditional expectations are positive maps of
\( * \)-algebras. As an illustration, consider the case that \( D\subset B \)
consists of scalar multiples of \( 1 \), and \( F:B\to D \) is such that \( \tau =F\circ E \)
is a \emph{trace} on \( M \), satisfying \( \tau (xy)=\tau (yx) \) for all
\( x,y\in M \).

Recall that \( X \) is called a \( B \)-semicircular variable if its cumulant
series is given by
\begin{equation}
k^{X}_{B;\underbrace{{1,\dots ,1}}_{p}}(b_{1},\dots ,b_{p-1})
=\kk^B(X,b_1X,\dots,b_{p-1}X)=\delta _{p,2}\eta (b_{1})
\end{equation}
for some map \( \eta :B\to B \). It is easily seen that if \( X \) is \( B \)-semicircular,
then \( \eta (b)=E(XbX) \).

\begin{thm}
\label{thrm:characteriRestrictedSemicircular}Let \( (M,E:M\to B) \) be a \( B \)-probability
space, such that \( M \) and \( B \) are \( C^{*} \)-algebras. Let \( F:B\to \mathbb {C}=D\subset B \)
be a faithful state. Assume that \( \tau =F\circ E \) is a faithful trace on
\( M \). Let \( X \) be a \( B \)-semicircular variable in \( M \). Then
the distribution of \( X \) with respect to \( \tau  \) is the semicircle
law iff \( E(X^{2})\in \mathbb {C} \).
\end{thm}
\begin{proof}
If \( E(X^{2})=k_{B;11}^{X}(1)\in \mathbb {C} \), it follows that the \( B \)-valued
cumulants of \( X \), restricted to \( D=\mathbb {C} \) are valued in \( D \).
Hence by Theorem \ref{thrm:restricRtransformDvalued}, the \( D \)-valued cumulant
series of \( X \) are the same as the restriction of the \( B \)-valued cumulant
series; hence the only scalar-valued cumulant of \( X \) which is nonzero is
the second cumulant \( k_{11}^{X} \), so that the distribution of \( X \)
is the semicircle law.

Conversely, assume that the distribution of \( X \) is the semicircle law.
Let \( \eta (b)=\kk^B(X,bX)=E(XbX) \), \( b\in B \). Then we have
\begin{align*}
2\tau (\eta (1))^{2} &= 2\tau (X^{2})^{2}\\
 &= \tau (X^{4})\\
 &= F\circ E(X^{4})
\end{align*}
Using the $B$-valued moment-cumulant formula for $X$, we get that
\begin{align*}
E(X^4)&=\mu^B(X,X,X,X)\\
&=\kk^B(X\kk^B(X,X),X)+\kk^B(X,X)\kk^B(X,X)\\
&=\eta(\eta(1))+\eta(1)\eta(1),
\end{align*}
and we can continue the above calculation as follows:
\begin{align*}
2\tau (\eta (1))^{2} 
 &= F(\eta (\eta (1)))+F(\eta (1)\eta (1))\\
 &= \tau (X\eta (1)X)+\tau (\eta (1)^{2})\\
 &= \tau (\eta (1)XX)+\tau (\eta (1)^{2})\\
 &= 2\tau (\eta (1)^{2}),
\end{align*}
so that
\[
\tau (\eta (1))=\tau (\eta (1)^{2})^{1/2}.\]
By the Cauchy-Schwartz inequality, we have that if \( \eta (1)\notin \mathbb {C} \),
\[
\tau (\eta (1))=\langle \eta (1),1\rangle <\Vert \eta (1)\Vert _{2}\cdot \Vert 1\Vert _{2}=\tau (\eta (1)^{2})^{1/2},\]
which is a contradiction. Hence \( \eta (1)\in \mathbb {C} \).
\end{proof}
We mention a corollary, which is of interest to random matrix theory. Let \( \sigma (x,y)=\sigma (y,x) \)
be a non-negative function on \( [0,1]^{2} \), having at most a finite number
of discontinuities in each vertical line. Let \( G(n) \) be an \( n\times n \)
random matrix with entries \( g_{ij} \), so that \( \{g_{ij}:i\leq j\} \)
are independent complex Gaussian random variables, \( g_{ij}=\overline{g_{ji}} \),
the expectation \( E(g_{ij})=0 \) and the variance \( E(|g_{ij}|^{2})=\frac{1}{n}\sigma (\frac{i}{n},\frac{j}{n}) \).
The matrices \( G(n) \) are called Gaussian Random Band Matrices. Let \( \mu _{n} \)
be the expected eigenvalue distribution of \( G(n) \), i.e.,
\[
\mu _{n}([a,b])=\frac{1}{n}\times \textrm{expected number of eigenvalues of }G(n)\, \textrm{in }[a,b].\]

\begin{cor}
The eigenvalue distribution measures \( \mu _{n} \) of the Gaussian Random
Band Matrices \( G(n) \) converge weakly to the semicircle law iff \( \int _{0}^{1}\sigma (x,y)dy \)
is a.e. a constant, independent of \( x \).
\end{cor}
The proof of this relies on a result from \cite{shlyakht:bandmatrix}, showing
that \( G(n) \) has limit eigenvalue distribution \( \mu  \), given as follows.
Let \( X \) be the \( L^{\infty }[0,1] \)-semicircular variable in an \( L^{\infty }[0,1] \)-probability
space \( (M,E:M\to L^{\infty }[0,1]) \), so that \( E(XfX)(x)=\int _{0}^{1}f(y)\sigma (x,y)dy \).
Let \( F:L^{\infty }[0,1]\to \mathbb {C} \) denote the linear functional \( F(f)=\int _{0}^{1}f(x)dx \),
and denote by \( \tau  \) the trace \( F\circ E \) on \( W^{*}(X,L^{\infty }[0,1]) \).
Then \( \mu  \) is the scalar-valued distribution of \( X \) with respect
to \( \tau  \), i.e.,
\[
\int t^{k}d\mu (t)=\tau (X^{k}).\]
It remains to apply Theorem \ref{thrm:characteriRestrictedSemicircular}, to
conclude that \( \mu  \) is a semicircle law iff \( E(X^{2})\in \mathbb {C} \),
i.e., \( \int _{0}^{1}\sigma (x,y)dy \) is a constant function of \( x \).

\subsection{A characterization of freeness.}

We are now ready to state the first main result of this note. The following theorem
was earlier proved for \( B \)-valued semicircular variables in \cite{shlyakht:amalg},
and found many uses in operator algebra theory.

\begin{thm}   \label{theorem:main}
Let \( X_{1},\dots ,X_{n}\in M \).   
Assume that $F:B\to D$ satisfies
the faithfulness condition that if $b_1\in B$ and 
if $F(b_1 b_2)=0$ for all $b_2\in B$, then 
$b_1 = 0$.
Then \( X_{1},\dots ,X_{n} \) are free
from \( B \) with amalgamation over 
\( D \) iff their \( B \)-valued cumulant
series satisfies
\begin{equation}
\label{eqn:RtransformFactorizationCondition}
k^{X_{1},\dots ,X_{n}}_{B;i_{1},\dots ,i_{k}}(b_{1},\dots ,b_{k-1})=
F\bigl(k^{X_{1},\dots ,X_{n}}_{B;i_{1},\dots ,i_{k}}(F(b_{1}),\dots ,F(b_{k-1}))
\bigr).
\end{equation}
for all \( b_{1},\ldots ,b_{k-1}\in B \). In short, \( k=F\circ k\circ F \).
(Here the cumulant series are computed in the \( B \)-probability space \( (M,E:M\to B) \)).
Equivalently,
\begin{equation}
\label{equiv}
k^{X_{1},\dots ,X_{n}}_{B;i_{1},\dots ,i_{k}}
(b_{1},\dots ,b_{k-1})=k^{X_{1},\dots ,X_{n}}_{D;i_{1},\dots ,i_{k}}
(F(b_{1}),\dots ,F(b_{k-1})).
\end{equation}
\end{thm}
We note that in the case that $M$ is a $C^*$-probability space,
the faithfulness assumption above is exactly the condition that the
GNS representation of $B$ with respect to the conditional expectation
$F$ is faithful.

Note also that the two characterizations in terms of cumulants 
appearing in the above theorem are
of a different nature: Eq. (\ref{eqn:RtransformFactorizationCondition}) 
is a condition ("factorization") on the
$B$-valued cumulants, whereas Eq. (\ref{equiv}) is a statement about the relation
between $B$-valued and $D$-valued cumulants. The equivalence of these two
formulations is easily seen with the help of Theorem \ref{restriction}.

Since \( X_{1},\dots ,X_{n} \) are free from \( B \) with amalgamation over
\( D \) iff the algebra \( N \) generated by \( X_{1},\dots ,X_{n} \) and
\( D \) is free from \( B \) over \( D \), the theorem above can be equivalently
stated as

\begin{thm}
Let \( N\subset M \) be a subalgebra of \( M \), containing \( D \). 
Assume that $F:B\to D$ satisfies
the faithfulness condition that if $b_1\in B$ and 
if $F(b_1 b_2)=0$ for all $b_2\in B$, then 
$b_1 = 0$.
Then
\( N \) is free from \( B \) over \( D \) iff
\begin{equation}
\kk^B(n_{1}b_{1},n_{2}b_{2},\dots ,n_{k})=
F\bigl(\kk^B(n_{1}F(b_{1}),n_{2}F(b_{2}),\dots ,n_{k})\bigr)
\end{equation}
for all \( k \) and all 
\( n_{1},\dots ,n_{k}\in N \), \( b_{1},\dots ,b_{k-1}\in B \). 
Equivalently,
\begin{equation}
\kk^B(n_{1}b_{1},n_{2}b_{2},\dots ,n_{k})=
\kk^D(n_{1}F(b_{1}),n_{2}F(b_{2}),\dots ,n_{k}).
\end{equation}

\end{thm}
\begin{proof}
We prove the theorem in the first formulation.

Assume that the condition (\ref{eqn:RtransformFactorizationCondition}) is satisfied
by the cumulant series of \( X_{1},\dots ,X_{n} \). Let \( Y_{1},\dots ,Y_{n} \)
be as in Section \ref{sec:canonical}. Since the freeness of \( X_{1},\dots ,X_{n} \)
from \( B \) with amalgamation over \( D \) is a condition on the \( B \)-moment
series of \( X_{1},\dots ,X_{n} \), and \( Y_{1},\dots ,Y_{n} \) have the
same \( B \)-moment series as \( X_{1},\dots ,X_{n} \), it is sufficient to
prove that \( Y_{1},\dots ,Y_{n}\in \mathcal{F} \) are free with amalgamation
over \( D \) from \( B \).

Since (\ref{eqn:RtransformFactorizationCondition}) is satisfied, \( \lambda _{j}^{0}=E_{B}(Y_{j})\in D \)
and hence \( Y_{1},\dots ,Y_{n} \) belong to the algebra \( \mathcal{L} \)
generated in \( \mathcal{F} \) by \( \lambda _{j}^{*} \) and \( \lambda _{p}^{q} \),
\( 1\leq j,p\leq n \), \( q\geq 1 \) and \( D \). Therefore, it is sufficient
to prove that \( \mathcal{L} \) is free from \( B \) with amalgamation over
\( D \).

Let \( w_{1},\dots ,w_{s}\in \mathcal{L} \), so that \( F\circ E_{B}(w_{j})=0 \),
and let \( b_{0},\dots ,b_{s}\in B \), so that \( F(b_{j})=0 \) (allowing
also \( b_{0} \) and/or \( b_{s} \) to be equal to \( 1 \)). We must prove
that
\[
F\circ E_{B}(b_{0}w_{1}b_{1}\cdots w_{s}b_{s})=0.\]
Note that the factorization condition (\ref{eqn:RtransformFactorizationCondition})
as well as the definition of the generators of \( \mathcal{L} \) (see \ref{relation:lambdas}
and \ref{relation:E}) imply that \( E_{B}|_{\mathcal{L}} \) has values in
\( D \). It follows that we may assume that \( E_{B}(w_{j})=0 \) (since \( F\circ E_{B}(w_{j})=E_{B}(w_{j})\in D \)).
By the definition of \( E_{B} \), its kernel is spanned by irreducible non-trivial
words in the generators \( \lambda _{j}^{*} \) and \( \lambda _{p}^{q} \).
Then
\[
W=b_{0}w_{1}\cdots w_{s}b_{s}\]
 is again a linear combination of words in the generators \( \lambda _{j}^{*} \)
and \( \lambda _{p}^{q} \). By linearity, we may reduce to the case that \( W \)
is a single word. If \( W \) is irreducible, it must be non-trivial (since
each \( w_{i} \) is non-trivial), hence \( E_{B}(W)=0 \), so that \( F\circ E_{B}(W)=0 \).
So assume that \( W \) is not irreducible. Since each \( w_{i} \) is irreducible,
this means that \( W \) contains a sub-word of the form
\begin{eqnarray*}
W & = & W_{1}\cdot d_{0}\lambda _{i_{1}}^{*}d_{1}\lambda ^{*}_{i_{2}}d_{2}\cdots \\
 &  & \quad d_{s}b_{1}d_{s+1}\lambda _{i_{s+2}}^{*}\cdots b_{j}\cdots \lambda _{i_{k}}^{*}d_{k}b_{r}d_{k+1}\lambda _{i_{k+1}}^{*}\\
 &  & \qquad \cdots d_{p}\lambda _{i_{p}}^{*}d_{p+1}\lambda _{j}^{p}\cdot W_{2}.
\end{eqnarray*}
Using the relation (\ref{relation:lambdas}) and the factorization condition
(\ref{eqn:RtransformFactorizationCondition}), we get that
\[
W=W_{1}d_{0}k_{i_{1},\dots ,i_{p},j}(d_{1},\dots ,F(d_{k}b_{r}d_{k+1}),\dots ,d_{p+1})W_{2}=0,\]
since \( F(d_{k}b_{r}d_{k+1})=d_{k}F(b_{r})d_{k+1}=0 \). Thus in any case,
\( F\circ E_{B}(W)=0 \).

We have therefore seen that the factorization condition implies freeness with
amalgamation.

To prove the other implication, assume that \( X_{1},\dots ,X_{n} \) are free
with amalgamation over \( D \) from \( B \). Let \( Z_{1},\dots ,Z_{n} \)
be \( B \)-valued random variables, so that
\[
k_{B;i_{1},\dots ,i_{k}}^{Z_{1},\dots ,Z_{n}}(b_{1},\dots ,b_{k-1})=
F\bigl(k^{X_{1},\dots ,X_{n}}_{D,i_{1},\dots ,i_{k}}(F(b_{1}),\dots ,F(b_{k-1}))
\bigr).\]
where \( k_{D} \) denote \( D \)-valued cumulants. (Note that the first occurrence
of \( F \) is actually redundant, as 
\( k^{X_{1},\dots ,X_{n}}_{D,i_{1},\dots ,i_{k}}(F(b_{1}),\dots ,F(b_{k-1}))\in D \)).
Then by the first part of the proof, \( Z_{1},\dots ,Z_{n} \) are free from
\( B \) with amalgamation over \( D \). Moreover, by Theorem \ref{thrm:restricRtransformDvalued},
the \( D \)-valued distributions of \( Z_{1},\dots ,Z_{n} \) and \( X_{1},\dots ,X_{n} \)
are the same. By assumption, \( X_{1},\dots ,X_{n} \) are free with amalgamation
over \( D \) from \( B \). This freeness, together with the \( D \)-valued
distribution of \( X_{1},\dots ,X_{n} \), determines their \( B \)-valued
distribution.  Indeed, the freeness assumptions determine
$$ F\circ E  (b' b_0 X_{i_1} b_1 X_{i_2} \cdots b_{n-1} X_{i_n} b_n),\quad b',b_i\in B$$
which in view of the assumptions on $F$ determines $$
 E  (b_0 X_{i_1} b_1 X_{i_2} \cdots b_{n-1} X_{i_n} b_n),\quad b_i\in B.$$
 It follows that the \( B \)-valued distributions of \( X_{1},\dots ,X_{n} \)
and \( Z_{1},\dots ,Z_{n} \) coincide. Hence the \( B \)-valued cumulants
of \( X_{1},\dots ,X_{n} \) satisfy (\ref{eqn:RtransformFactorizationCondition}). 
\end{proof}

As an application, we have the following proposition. 

\begin{prop}
\label{application}
Let \( N\subset M \) be a subalgebra. Let \( D\subset C\subset B\subset M \)
be subalgebras, and \( E_{C}:M\to C \), \( E_{B}:M\to B \), \( E_{D}:M\to D \)
be conditional expectations, so that 
\( E_{D}=E_{D}\circ E_{C} \), \( E_{C}=E_{C}\circ E_{B} \).  
Assume that $E_C\vert_B:B\to C$, $E_D\vert_C:C\to D$ 
and $E_D\vert_B:B\to D$ satisfy
the faithfulness assumptions of Theorem~\ref{theorem:main}.
Assume that \( N \) is free from \( C \) with amalgamation over \( D \),
and also free from \( B \) with amalgamation over \( C \). Then \( N \) is
free from \( B \) with amalgamation over \( D \).
\end{prop}

\begin{proof}
Since \( N \) is free from \( B \) with amalgamation over \( C \), we have
that for all \( n_{j}\in N \) and \( b_{j}\in B \),
\[
\kk^B(n_{1}b_{1},n_{2}b_{2},\dots ,n_{k})=\kk^C(
n_{1}E_{C}(b_{1}),\dots ,n_{k}).\]
Since \( N \) is free from \( C \) with amalgamation over \( D \), we get
similarly that for all \( c_{j}\in C \),
\[
\kk^C(n_{1}c_{1},n_{2}c_{2},\dots ,n_{k})=\kk^D
(n_{1}E_{D}(c_{1}),\dots ,n_{k}).\]
Applying this with \( c_{j}=E_{C}(b_{j}) \) and combining with the previous
equation gives
\[
\kk^B(n_{1}b_{1},n_{2}b_{2},\dots ,n_{k})=\kk^D(
n_{1}E_{D}(b_{1}),\dots ,n_{k}),\]
since \( E_{D}=E_{D}\circ E_{C} \). Hence \( N \) is free from \( B \) with
amalgamation over \( D \).
\end{proof}

In general, for operator-valued random variables $X$ and $Y$, freeness 
over $B$ and
freeness over $D$ are not implying
each other in a straightforward manner. What can be said in general is
that freeness of $B\langle X\rangle$ from Y over $D$ implies freeness
of $X$ from $Y$ over $B$ (compare \cite[Lemma 2.6]{shlyakht:cpentropy}).
It is therefore interesting
to note that the above proposition provides us with a situation where
we also get the reverse implication.

\begin{cor}
\label{cor:free}
Let $(M,\tau)$ be a tracial $W^*$-probability space (i.e., $M$ is a von
Neumann algebra and $\tau:M\to\CC$ a faithful and normal trace).
Let $D\subset B\subset M$ be two von Neumann subalgebras and consider two
selfadjoint random variables $X,Y\in M$. Assume that $Y$ is free from $B$
over $D$. Then the following two statements are equivalent:
\begin{enumerate}
\item
$X$ and $Y$ are free with amalgamation over $B$.
\item
$B\langle X\rangle$ and $Y$ are free with amalgamation
over $D$.
\end{enumerate}
\end{cor}

Of course, all freeness statements above are with respect to the unique
trace-preserving conditional expectations.

\begin{proof}
That freeness over $D$ implies freeness over $B$ is the statement of
Lemma 2.6 of \cite{shlyakht:cpentropy}. The reverse implication is a direct
application of our Prop.~\ref{application}. One should note that the faithfulness
assumptions are automatically fulfilled in our frame, where all conditional expectations
are compatible with the faithful state $\tau$.
\end{proof} 

Let us also mention the following corollary of Prop. \ref{application}.

\begin{cor}
Let \( D\subset C\subset B \)
be \( C^{*} \)-algebras, \( E^{N}_{D}:N\to D \), \( E^{C}_{D}:C\to D \) and
\( E^{B}_{C}:B\to C \) be conditional expectations, 
having faithful GNS representations.
Consider the reduced free product
\[
(((N,E^{N}_{D})*_{D}(C,E^{C}_{D})),E^{N}_{D}*\textrm{id})*_{C}(B,E^{B}_{C}),\]
where \( E^{N}_{D}*\textrm{id} \) denotes the canonical conditional expectation
from the free product \( (N,E^{N}_{D})*_{D}(C,E^{C}_{D}) \) onto \( C \).
Then
\begin{equation}
(((N,E^{N}_{D})*_{D}(C,E^{C}_{D})),E^{N}_{D}*\textrm{id})*_{C}(B,E^{B}_{C}) 
\cong (N,E^{N}_{D})*_{D}(B,E^{C}_{D}\circ E^{B}_{C}).  
\end{equation}
\end{cor}

\begin{proof}
To see this, it is sufficient to prove that 
\( N\subset (((N,E^{N}_{D})*_{D}(C,E^{C}_{D})),
E^{N}_{D}*\textrm{id})*_{C}(B,E^{B}_{C}) \)
is free from \( B \) with amalgamation over 
\( D \), since both \( (((N,E^{N}_{D})*_{D}(C,E^{C}_{D})),
E^{N}_{D}*\textrm{id})*_{C}(B,E^{B}_{C}) \)
and \( (N,E^{N}_{D})*_{D}(B,E^{C}_{D}\circ E^{B}_{C}) \) are generated by \( N \)
and \( B \) as \( C^{*} \)-algebras. But \( N \) is free from \( B \) over
\( C \), and from \( C \) over \( B \), by construction. Hence by the proposition
above, \( N \) is free from \( B \) over \( D \).
\end{proof}

\section{Relative Fisher information and conjugate variables with respect to
a subalgebra}
\subsection{Basic definitions}
 From now on we will work in a \emph{tracial $W^*$-probability space}
$(M,\tau)$, i.e. $M$ is a von Neumann algebra and $\tau:M\to\CC$ a faithful
and normal trace. For given subsets $S_1,\dots,S_p\subset M$, we will 
denote by $L^2(S_1,\dots,S_p)$ the closure of the von Neumann subalgebra
generated by all $S_i$ ($i=1,\dots,p$) in the $L^2$-norm of $\tau$,
i.e. $\Vert a\Vert_{L^2}^2=\tau(aa^*)$. (Usually, the sets $S_i$ will
be either subalgebras or consist of given random variables.)

Let $B\subset M$ 
be a unital $*$-subalgebra and consider selfadjoint random variables $X_1,\dots,X_n\in M$.
Recall \cite{V5} that 
the \emph{conjugate variables of $X_1,\dots,X_n$ with respect to $B$},
$J_1=J_1(X_1,\dots,X_n:B),\dots,
J_n=J_n(X_1,\dots,X_n:B)$,
are determined by the requirements that they belong to $L^2(X_1,\dots,X_n,B)$ and
that they fulfill the following system of equations:
\begin{multline}
\label{conjugate}
\tau(J_ib_1X_{i_1}b_2X_{i_2}b_3\cdots b_{m}X_{i_m}b_{m+1})=\\ \sum_{r=1}^m
\delta_{ii_r}\cdot
\tau\bigl(b_1X_{i_1}b_2\cdots b_{r-1}X_{i_{r-1}}b_{r}\bigr)
\cdot\tau\bigl(b_{r+1} X_{i_{r+1}}b_{r+2}\cdots
 b_{m}X_{i_m}b_{m+1}\bigr)
\end{multline}
for all $m\geq 0$, all $b_1,\dots,b_{m+1}\in B$, and all $1\leq i_1,\dots
,i_m\leq n$. The equations for $m=0$ are to be understood as
\begin{equation}
\tau(J_ib)=0
\end{equation}
for $i=1,\dots,n$ and all $b\in B$.

Furthermore, if a system of conjugate variables $J_1,\dots,J_n$
for $X_1,\dots,X_n$ with respect
to $B$ exists (in which case it is unique and satisfies $J_i=J_i^*$ for all 
$i=1,\dots,n$), then Voiculescu defined
\begin{equation}
\Phi^*(X_1,\dots,X_n:B):=\sum_{i=1}^n\Vert J_i\Vert_{L^2(\ff)}^2=
\sum_{i=1}^n \ff(J_i^2)
\end{equation}
and called it the \emph{relative free information with respect to $B$}.
If no system of conjugate variables exist, then he put $\Phi^*(X_1,\dots,X_n:B):=\infty$.

\subsection{Reformulation in terms of cumulants}
It is quite easy to check that the above equations can also be rewritten in
terms of scalar-valued cumulants $\kk=\kk^\CC$ in the equivalent form:
\begin{align}
\kk(J_i)&=0\\
\kk(J_i,a)&=\delta_{a,X_i}\\
\kk(J_i,a_1,\dots,a_m)&=0
\end{align}
for all $m\geq 2$, all $i=1,\dots,n$ and all $a,a_1,\dots,a_m\in\{X_1,\dots,X_n\}
\cup B$.

If we are in the context of a $B$-valued probability space $(M,E:M\to B)$ which
is compatible with the scalar-valued probability structure, i.e. $\tau\circ E=\tau$, then
the definition of the 
conjugate variables $J_i(X_1,\dots,X_n:B)$ has an operator-valued flavour and
it seems conceivable that there should also exist a nice description of the determining
equations in terms of $B$-valued cumulants. This is indeed the case, namely, 
as we will see below, we have:
\begin{align}
\kk^B(J_i)&=0\\
\kk^B(J_i,bX_j)&=\delta_{ij}\tau(b)\\
\kk^B(J_i,b_1X_{i_1},b_2X_{i_2},\dots,b_mX_{i_m})&=0
\end{align}
for all $m\geq 2$, all $i=1,\dots,n$, all $b,b_1,\dots,b_m\in B$, and all 
$1\leq j,i_1,\dots,i_m\leq n$.
Note the change of the role of the elements from $B$:
from arguments in the scalar-valued cumulants to (non-commuting) constants
in the $B$-valued cumulants.

These descriptions of the conjugate variables with respect to $B$ in terms
of the $\CC$-valued and in terms of the $B$-valued distribution are the extreme
cases of the following more general description in terms of the $D$-valued
distribution for any intermediate subalgebra $\CC\subset D\subset B$.

\begin{thm}
\label{conjugate-D}
Let $(M,\tau)$ be a tracial $W^*$-probability space
and $B\subset M$ a unital $*$-subalgebra. Consider 
selfadjoint random variables $X_1,\dots,X_n\in M$ and their conjugate variables
$J_i:=J_i(X_1,\dots,X_n:B)$ with respect to $B$.
Let $D\subset B$ be a unital subalgebra of $B$ with conditional
expectation $E:M\to D$ such that $\tau\circ E=\tau$.
Then the defining equations (\ref{conjugate}) for the $J_i$ are equivalent to
the following system of equations in terms of $D$-valued cumulants:
\begin{align} 
\kk^D(J_i)&=0 \label{eins}\\
\kk^D(J_i,da)&=\delta_{aX_i}\tau(d)\\
\kk^D(J_i,d_1a_1,\dots,d_ma_m)&=0\label{drei}
\end{align}
for all $m\geq 2$, all $i=1,\dots,n$, all $d,d_1,\dots,d_m\in D$,
and all $a,a_1,\dots,a_m\in\{X_1,\dots,X_n\}
\cup B$.
\end{thm}

Note that the traciality of $\tau$ implies that we have for
all $r\geq 2$ and all $m_1,\dots,m_r\in M$
\begin{equation}
\label{traciality}
\tau\bigl(\kk^D(m_1,\dots,m_r)\bigr)=\tau\bigl(\kk^D(m_2,\dots,m_r,m_1)\bigr).
\end{equation} 
Together with the faithfulness of $\tau$ this yields that we have the same
kind of 
formulas as in Theorem \ref{conjugate-D} also in the cases where $J_i$ is not
the first argument of a cumulant $\kk^D$, but appears at an arbitrary position.

\begin{proof}
It is easy to see that, by the moment-cumulant formula (\ref{moment-cumulant}),
the Equations (\ref{eins}) - (\ref{drei}) are equivalent to
$$E(J_id_1a_1\cdots d_ma_m)=\sum_{r=1}^m\delta_{a_rX_i}\cdot \tau(d_1a_1\cdots d_{r})\cdot
E(d_{r+1}a_{r+1}\cdots d_ma_m)$$
for all $m\geq 0$, all $d_1,\dots,d_m\in D$ and all $a_1,\dots,a_m\in\{X_1,
\dots,X_n\}\cup B$.
Since $\tau$ is faithful and $\tau\circ E=\tau$, this is equivalent to
$$\tau(J_id_1a_1\cdots d_ma_md_{m+1})=\sum_{r=1}^m\delta_{a_rX_i}\cdot
\tau(d_1a_1\cdots d_{r})\cdot
\tau(d_{r+1}a_{r+1}\cdots d_ma_md_{m+1})$$
for all $m\geq 0$, all $d_1,\dots,d_m,d_{m+1}\in D$ and all $a_1,\dots,a_m\in\{X_1,
\dots,X_n\}\cup B$.
But this is clearly the same as (\ref{conjugate}).
\end{proof}

\subsection{Relation between $\Phi^*(X:D)$ and $\Phi^*(X:B)$}
Our aim will be to investigate, for given random variables $X_1,\dots,X_n\in M$,
the relation between their Fisher informations with respect to two different
subalgebras. If we have $D\subset B\subset M$, then the following facts are known
from \cite{V5,V6,shlyakht:cpentropy}:
\begin{enumerate}
\item
We always have: $\Phi^*(X_1,\dots,X_n:D)\leq \Phi^*(X_1,\dots,X_n:B)$.
\item
If $\{X_1,\dots,X_n\}$ is free from $B$ with amalgamation over $D$, then we
have equality: $\Phi^*(X_1,\dots,X_n:D)= \Phi^*(X_1,\dots,X_n:B)$.
\item
If $D=\CC$, then the reverse implication of (ii) holds: If $\Phi^*(X_1,\dots,X_n:B)=
\Phi^*(X_1,\dots,X_n)<\infty$, then we have that $\{X_1,\dots,X_n\}$ is free from
$B$ (i.e. free from $B$ with amalgamation over $\CC$).
\end{enumerate}
The question which we want to address is whether the statement (iii) also holds
for more general $D$, i.e. is it true that
$\Phi^*(X_1,\dots,X_n:B)=\Phi^*(X_1,\dots,X_n:D)$ implies that $\{X_1,\dots,X_n\}$
is free from $B$ with amalgamation over $D$.
We will be able to show that this is true for finite-dimensional $D$. 

The techniques for proving this are operator-valued
generalizations of Voiculescu's ideas \cite{V6} for
dealing with the special case $D=\CC$. The main conceptual ingredient will be
an operator-valued version of the liberation gradient, which we present in the
next section.

\section{Operator-valued liberation gradient}
\subsection{Definition and basic properties}
In \cite{V6}, Voiculescu introduced the notion $\ff^*(A_1,A_2)$ 
of liberation Fisher information,
which is a measure for how far two $*$-subalgebras $A_1$ and $A_2$ 
in a tracial $W^*$-probability space are away from being free. As in the case
of the free Fisher information this liberation Fisher information is given by
the square of the $L^2$-norm of a special vector in $L^2(A_1,A_2)$, namely of
the so-called liberation gradient $j(A_1:A_2)$. The defining
property for this liberation gradient $j$ is in terms of a canonical
derivation $\delta$.
The definition of Voiculescu is recovered as the case $B=\CC$ and $E=\tau$ of
our following operator-valued generalizations.

\begin{defn}
\label{def:gradient}
Let $(M,\tau)$ be a tracial $W^*$-probability space, and let $E:M\to B$ be
a conditional expectation with $\tau\circ E=\tau$. Consider two subalgebras $A_1,A_2$,
both of them containing $B$ as a subalgebra, 
$B\subset A_1,A_2 \subset M$, which are algebraically free modulo $B$ (i.e.,
the canonical homomorphism $A_1*_BA_2\to A_1\vee A_2$ has trivial kernel). 
\\
1) Denote by $A=A_1\vee A_2$ 
the algebra generated
by $A_1$ and $A_2$. 
We define
\begin{equation}
\delta_{A_1:A_2}:A\to A\otimes_B A
\end{equation}
to be the derivation into the $A$-bimodule $A\otimes_B A$ which is
determined by
\begin{equation}
\delta_{A_1:A_2}(a)=\begin{cases}
a\otimes 1-1\otimes a, &\text{if $a\in A_1$}\\
0,&\text{if $a\in A_2$}
\end{cases}
\end{equation}
\\
2) We define the \emph{$B$-valued liberation gradient}
$j:=j_B(A_1:A_2)$ by the requirements that it is a vector in $L^2(A_1,A_2)$
and that we have for all $a\in A_1\vee A_2$
\begin{equation}
\label{gradient}
E(ja)=E\otimes E(\delta_{A_1:A_2}(a)).
\end{equation}
3) We define the \emph{$B$-valued liberation Fisher information} by
\begin{equation}
\ff^*_B(A_1:A_2):=
\begin{cases}
\Vert j_B(A_1:A_2)\Vert_{L^2(\tau)}^2,&\text{if $j_B(A_1:A_2)$ exists}\\
\infty,&\text{otherwise}
\end{cases} 
\end{equation}
\end{defn}

\begin{rems}
\label{rem:delta}
1) It is easy to see that the relations $E(ja)=E\otimes E(\delta_B (a))$
have the following explicit form:
\begin{multline}
\label{gradient2}
E(jc_1\tilde c_1c_2\tilde c_2\dots c_m\tilde c_m)=\sum_{r=1}^m \bigl( E(c_1
\tilde c_1\dots
\tilde c_{r-1}c_r)E(\tilde c_rc_{r+1}\tilde c_{r+1}\cdots c_m\tilde c_m)\\-
 E(c_1\tilde c_1\dots
c_{r-1}\tilde c_{r-1})E(c_r\tilde c_r\cdots c_m\tilde c_m)\bigr)
\end{multline}
for all $m\geq 0$ and all $c_1,\dots,c_m\in A_1$, $\tilde c_1,\dots,\tilde c_m\in A_2$.
(The case $m=0$ has then to be interpreted as $E(j)=0$.)
\\
2) Note that if a liberation gradient exists, then it is uniquely determined by the above
requirements. 
\\
3) Since, for $b\in B=A_1\cap A_2$, we have 
$$0=\delta_{A_1:A_2}(b)=b\otimes 1-1\otimes b,$$
it is clear that $\delta_{A_1:A_2}$ has to take values in $A\otimes_B A$, and not
just in $A\otimes A$. 
\end{rems}

It will be quite instructive to characterize the liberation gradient in terms
of $B$-valued cumulants. 

\begin{thm}
\label{char-gradient}
The equations (\ref{gradient}) in the definition of the operator-valued liberation
gradient are equivalent to the following system of equations in terms of
$B$-valued cumulants: 
\begin{equation}
\kk^B(j)=0
\end{equation}
and,
for all $m\geq 1$ and all 
$a_1,\dots,a_m\in A_1\cup A_2$,
\begin{equation}
\label{eq-gradient}
\kk^B(j,a_1,\dots,a_m)=\begin{cases}
0,&\text{if $a_1,a_m\in A_1$ or if $a_1,a_m\in A_2$}\\
-\kk^B(a_1,\dots,a_m),&\text{if $a_1\in A_1$ and $a_m\in A_2$}\\
+\kk^B(a_1,\dots,a_m),&\text{if $a_1\in A_2$ and $a_m\in A_1$}
\end{cases}
\end{equation}
\end{thm}

\begin{proof}
Let us denote in the following by 
$\kk^B\otimes\kk^B$ the family of multi-linear mappings with one
argument from $A\otimes_B A$ and all other arguments from $A$, which
are determined as follows:
\begin{multline}
\kk^B\otimes \kk^B(a_1,\dots,a_{i-1},a\otimes \hat a,a_{i+1},\dots,a_m):=\\
\kk^B(a_1,\dots,a_{i-1},a)\cdot \kk^B(\hat a,a_{i+1},\dots,a_m).
\end{multline}
for all $m\geq 1$, $1\leq i\leq m$, $a_1,\dots,a,\hat a,\dots,a_m\in A$.
Then, by using the moment-cumulant formula (\ref{moment-cumulant}), one
can check that the system of equations
$E(ja)=E\otimes E(\delta (a))$
is equivalent to the system of equations
\begin{equation}
\label{prop-derivation}
\kk^B(j,a_1,\dots,a_m)=
\sum_{i=1}^m \kk^B\otimes\kk^B(a_1,\dots,a_{i-1},\delta a_i,a_{i+1},\dots,a_m)
\end{equation}
for all $m\geq 0$ and all $a_1,\dots,a_m\in A_1\vee A_2$.

If the $a_i$ appearing in Eq. (\ref{prop-derivation}) are from $A_1\cup A_2$, then this
equation reduces drastically due to the following observations: If $i\not=1$
and $i\not=m$, then we have for $a_i\in A_1$ 
\begin{align*}
\kk^B(\dots,\delta a_i,\dots)&=\kk^B(\dots,a_i\otimes 1-1\otimes a_i,\dots)\\
&=\kk^B(\dots,a_i)\cdot \kk^B(1,\dots) - \kk^B(\dots,1)\cdot \kk^B(a_i,\dots),
\end{align*}
which is equal to zero, because cumulants of length greater than 1 where
one entry is equal to 1 must vanish. Since for $a_i\in A_2$ the term $\delta a_i$
vanishes always, we have in any case that the terms with $1<i<m$ vanish.
On the other hand, if $m>1$, we have for the two remaining terms 
$$\kk^B(\delta a_1,a_2,\dots,a_m)=\begin{cases}
-\kk^B(a_1,a_2,\dots,a_m),&\text{if $a_1\in A_1$}\\
0,&\text{if $a_1\in A_2$}
\end{cases}$$
and 
$$\kk^B(a_1,a_2,\dots,\delta a_m)=\begin{cases}
\kk^B(a_1,a_2,\dots,a_m),&\text{if $a_m\in A_1$}\\
0,&\text{if $a_m\in A_2$}
\end{cases}.$$
If $m=1$, then we have $\kk^B(\delta a)=0$ for all $a\in A_1\cup A_2$.

Putting all these observations together gives the assertion.
\end{proof}

This theorem yields directly, by the fact that freeness is equivalent
to the vanishing of mixed cumulants
(see Theorem \ref{mixed}), the following
fundamental characterization of freeness with amalgamation in terms
of the operator-valued liberation gradient.

\begin{cor}
\label{cor:amalg}
Let $(M,\tau)$ be a tracial $W^*$-probability space, and let $E:M\to B$ be
a conditional expectation with $\tau\circ E=\tau$. 
Consider two subalgebras $A_1,A_2$,
both of them containing $B$ as a subalgebra, 
$B\subset A_1,A_2 \subset M$, which are algebraically free modulo $B$.
Then the following two statements are equivalent:
\begin{enumerate}
\item
$A_1$ and $A_2$ are free with amalgamation over $B$.
\item
$j_B(A_1:A_2)=0$
\end{enumerate} 
\end{cor}

\begin{rems}
1) The characterization of $j$ in terms of cumulants shows quite clearly
that $j$ contains the information about the mutual position of $A_1$ and $A_2$,
but that it ignores the information about the internal structure of $A_1$ and of
$A_2$.
\\
2) If $j_B(A_1:A_2)$ exists, then so does $j_B(A_2:A_1)$, and we have
the equality
\begin{equation}
j_B(A_1:A_2)=-j_B(A_2:A_1).
\end{equation}
This could be easily checked from the definition, but is strikingly clear
from Eq. (\ref{eq-gradient}).
\\
3) It is clear that $j_B(A_1:A_2)$ fulfills also the defining relations
for $j_B(\tilde A_1,\tilde A_2)$ if $\tilde A_1\subset A_1$ and $\tilde A_2
\subset A_2$. Thus we get $j_B(\tilde A_1:\tilde A_2)$ in this situation
by projecting $j_B(A_1:A_2)$ onto $L^2(\tilde A_1,\tilde A_2)$. 
\\
4) Instead of changing $A_1$ and $A_2$ it will be more relevant for our
questions how $j_B(A_1:A_2)$ behaves if we change the subalgebra $B$
(a more precise description of the framework for that question will be
given in Theorem \ref{relation}).
A first hint that there might be a relation between $j_D$ and $j_B$
is given by the following observation: 
Assume we have two unital subalgebras with $D\subset B$, such that
$E_D\circ E_B=E_D$. Then 
the defining
relations for $j:=j_D(A_1,A_2)$, $E_D(ja)=E_D\otimes E_D(\delta a)$, can also be
formulated equivalently with respect to $E_B$:
\begin{equation}
E_B(ja)=E_D\otimes E_B(\delta a).
\end{equation}
\end{rems}

\begin{prop}
\label{commutant}
If $j_B(A_1:A_2)$ exists then it belongs to the $\Vert\cdot\Vert_2$-closure
of the relative commutant $B'\cap M$. 
\end{prop}

\begin{proof}
We have to show that, for any $b\in B$, we have the equality
$$\tau(jbc_1\tilde c_1\cdots c_m\tilde c_m)=\tau(jc_1\tilde c_1\cdots c_m\tilde c_mb)$$
for all $m\geq 0$ and all $c_1,\dots,c_m\in A_1$, $\tilde c_1,\dots,\tilde c_m\in A_2$.
But this follows directly from applying $\tau$ to the equation (\ref{gradient2}).
\end{proof}

\begin{rems}
1) The last proposition shows that the definition of $j_B$ is in general
quite restrictive. In the case of infinite-dimensional $B$ the relative
commutant of $B$ might just consist of scalar multiples of the identity,
in which case we are only left with the dichotomy that either $j_B(A_1:A_2)=0$
(and thus $A_1$, $A_2$ are free with amalgamation over $B$) or they are so far
apart from being free over $B$ that no $j$ exists, i.e. the $B$-valued liberation
Fisher information is either zero or infinity.
\\
2) Prop.~\ref{commutant} 
suggests that, in the situation $D\subset B$, one might get $j_B$
from $j_D$ by projecting it onto the relative commutant of $B$. In
the case of finite dimensional $D$ we will make this rigorous in the next
section. 
\end{rems}

\subsection{Relation between $j$ and $J$ in the case of finite-dimensional $D$}
In this section we will make the additional assumption that $D$ is 
finite-dimensional. This has the effect that there exist explicit formulas for
the conditional expectations $E_D$ and $E_{D'}$ onto $D$ and onto the 
relative commutant $D'\cap M$, respectively.

We denote the group of unitary elements in $D$ by $U$. This is a compact group,
and we will denote integration with respect to its normalized Haar measure by 
"du".
Then there exists a positive, invertible and central element of $D$, which we
will denote by $c$ in the following, such that the following integral formulas hold
for all $m\in M$:
\begin{equation}
\label{integral-D}
E_D(m)=\dim(D)\cdot c^{-1}\cdot\int_U u\tau(u^*m)du
\end{equation}
\begin{equation}
\label{cond-D'}
E_{D'}(m)=\int_U umu^*du.
\end{equation}
The role of $c$ in the formula 
(\ref{integral-D}) is to correct the way how $\tau$ partitions
the unit between the minimal central projections of $D$. To be precise, suppose that
$$D\cong M_{n_1}(\CC)\oplus\dots\oplus M_{n_k}(\CC),$$
and let $p_1,\dots,p_k$ denote the 
minimal central projections of $D$. Then $c$ is given
as
$$c:=(n_1^2+\dots+n_k^2)\cdot\sum_{j=1}^k \frac{\tau(p_j)}{n_j^2}p_j.$$
The verification of the above integration formulas is straightforward.
They allow us to formulate and prove our main result about the connection between
the liberation gradient and the conjugate variable.
This will be a corollary of the following statement which clarifies
the relation between $j_\CC$ and $j_D$.

\begin{thm}
\label{relation}
Let $(M,\tau)$ be a tracial $W^*$-probability space
and $D\subset M$ a unital $*$-subalgebra. Consider 
two unital subalgebras $A_1,A_2\subset M$ with
$D\subset A_2$ (but not $D\subset A_1$).
Assume that $\dim(D)<\infty$.
If $j_\CC(A_1:A_2)$ exists, then $j_D(A_1\vee D:A_2)$
exists, too, and is given by
\begin{equation}
\label{lib-proj}
j_D(A_1\vee D:A_2))=E_{D'}\bigl( j_\CC(A_1:A_2)
\bigr)\cdot c^{-1}\cdot\dim(D).
\end{equation}
\end{thm}

\begin{proof}
Let us denote
$$j:=E_{D'}\bigl( j_\CC(A_1:A_2)
\bigr)\cdot c^{-1}\cdot\dim(D)=\int u j_\CC(A_1:A_2)u^*du
\cdot c^{-1}\cdot\dim(D).$$
Since this belongs to $L^2(D,A_1,A_2)$, it only remains to 
check the defining relations (\ref{gradient}) for $j$.
Let $a\in A_1\vee A_2$, then we have (with $E=E_D$)
\begin{align*}
E(ja)&=\int E(u j_\CC(A_1:A_2) u^* a)du 
\cdot c^{-1}\cdot\dim(D)\\
&=\int u\cdot E(j_\CC(A_1:A_2) u^* a)du\cdot c^{-1}\cdot\dim(D)\\
&=\int u \cdot\tau\otimes E(\delta_{A_1:A_2}(u^*a))du \cdot c^{-1}\cdot\dim(D)\\
&=\int u\cdot\tau\otimes E(u^* \delta_{A_1:A_2}(a))du \cdot c^{-1}\cdot\dim(D)\\
&=E\otimes E(\delta_{A_1\vee D:A_2}(a)),
\end{align*}
which yields the assertion.
\\
(One should note that $\delta_{A_1:A_2}(u^*)=0$, because $u^*\in A_2$, and that
on $A_1\vee A_2$ we can identify $\delta_{A_1:A_2}$ canonically with
$\delta_{A_1\vee D:A_2}$.)
\end{proof}

In the following, we will denote, for given random variables $X_1,\dots,X_n$
and a subalgebra $D$, by $D\langle X_1,\dots,X_n\rangle$ the algebra
generated by $D$ and $X_1,\dots,X_n$;
the elements of $D\langle X_1,\dots,X_n\rangle$
can be considered as non-commutative
polynomials in $X_1,\dots,X_n$ with coefficients from $D$.

\begin{cor}
\label{liberation-conjugate}
Let $(M,\tau)$ be a tracial $W^*$-probability space
and $B\subset M$ a unital $*$-subalgebra. Consider 
selfadjoint random variables $X_1,\dots,X_n\in M$ 
with $\Phi^*(X_1,\dots,X_n:B)<\infty$.
Let $D\subset B$ be a unital subalgebra of $B$ with conditional
expectation $E:M\to D$ such that $\tau\circ E=\tau$. Assume that
$\dim(D)<\infty$.  
Then the $D$-valued liberation gradient of the pair $D\langle X_1,\dots,X_n\rangle$ 
and $B$ exists and is given by
\begin{equation}
\label{lib-con}
j_D(D\langle X_1,\dots,X_n\rangle:B)=E_{D'}\bigl( \sum_{i=1}^n[J_i(X_1,\dots,X_n:B),X_i]
\bigr)\cdot c^{-1}\cdot\dim(D).
\end{equation}
\end{cor}

\begin{proof}
This follows from combining the above theorem with Voiculescu's 
formula \cite{V6},
$$j_\CC(\CC\langle X_1,\dots,X_n\rangle:B)=
\sum_{i=1}^n[J_i(X_1,\dots,X_n:B),X_i].$$
\end{proof}

\begin{cor}
\label{oben}
Let $(M,\tau)$ be a tracial $W^*$-probability space
and $B\subset M$ a unital $*$-subalgebra. Consider 
selfadjoint random variables $X_1,\dots,X_n\in M$ 
with $\Phi^*(X_1,\dots,X_n:B)<\infty$.
Let $D\subset B$ be a unital subalgebra of $B$ with conditional
expectation $E:M\to D$ such that $\tau\circ E=\tau$. Assume that
$\dim(D)<\infty$.  
Then the following two statements are equivalent:
\begin{enumerate}
\item 
$E_{D'}\bigl( \sum_{i=1}^n[J_i(X_1,\dots,X_n:B),X_i]
\bigr)=0$.
\item
$\{X_1,\dots,X_n\}$ is free from $B$ with amalgamation over $D$.
\end{enumerate}
\end{cor}

\begin{proof}
This is just a combination of Theorem \ref{liberation-conjugate}
with Corollary \ref{cor:amalg}.
\end{proof}

\begin{cor}
\label{unten}
Let $(M,\tau)$ be a tracial $W^*$-probability space. Let $D\subset M$ 
be a unital $*$-subalgebra with conditional expectation $E:M\to D$ such 
that $\tau\circ E=\tau$.
Assume that
$\dim(D)<\infty$. Consider 
selfadjoint random variables $X_1,\dots,X_n\in M$ and assume that
$\Phi^*(X_1,\dots,X_n:D)<\infty$. 
Then we have
\begin{equation}
E_{D'}\bigl( \sum_{i=1}^n[J_i(X_1,\dots,X_n:D),X_i]
\bigr)=0.
\end{equation}
\end{cor}

\begin{proof}
This follows from the fact that the freeness condition in
Corollary \ref{oben} is trivially  fulfilled for $B=D$.
\end{proof}

\section{Freeness from a subalgebra and equality of Fisher informations}

\subsection{Main results}
We are now ready to state our second main result of this note.

\begin{thm}
\label{main2a}
Let $(M,\tau)$ be a tracial $W^*$-probability space
and $A,B\subset M$ unital $*$-subalgebras which are algebraically
free. 
Let $D\subset B$ be a unital $*$-subalgebra of $B$ with conditional
expectation $E:M\to D$ such that $\tau\circ E=\tau$.   
Assume that $\dim D<\infty$.
If $j_\CC(A:B)$ exists,
then the following statements are
equivalent:
\begin{enumerate}
\item
$j_\CC(A:D)=j_\CC(A:B)$.
\item
$A$ is free from $B$ with amalgamation over $D$.
\end{enumerate}
\end{thm}

\begin{proof}
$(ii)\Rightarrow (i)$: This is Prop.~5.13 from \cite{V6}. 
\\
$(i)\Rightarrow (ii)$:
Theorem \ref{liberation-conjugate}
and Corollary \ref{unten} yield
\begin{align*}
j_D(A\vee D:B)&=
E_{D'}\bigl(j_\CC(A:B)
\bigr)\cdot c^{-1}\cdot\dim(D)\\
&=E_{D'}\bigl( j_\CC(A:D)
\bigr)\cdot c^{-1}\cdot\dim(D)\\
&=j_D(A\vee D:D)\\
&=0,
\end{align*}
where the vanishing of $j_D(A\vee D:D)$ follows from the fact that $A$ and $D$
are clearly free with amalgamation over $D$.
\\ 
This gives statement $(ii)$ by the fact that the vanishing of the liberation
gradient is equivalent to freeness (Corollary \ref{cor:amalg}).
\end{proof}

\begin{thm}
\label{main2}
Let $(M,\tau)$ be a tracial $W^*$-probability space
and $B\subset M$ a unital $*$-subalgebra. Consider 
selfadjoint random variables $X_1,\dots,X_n\in M$ 
with $\Phi^*(X_1,\dots,X_n:B)<\infty$.
Let $D\subset B$ be a unital $*$-subalgebra of $B$ with conditional
expectation $E:M\to D$ such that $\tau\circ E=\tau$.   
Assume that $\dim D<\infty$.
Then the following statements are
equivalent:
\begin{enumerate}
\item
$\Phi^*(X_1,\dots,X_n:D)=\Phi^*(X_1,\dots,X_n:B)$
\item
$\{X_1,\dots,X_n\}$ is free from $B$ with amalgamation over $D$.
\end{enumerate}
\end{thm}

\begin{proof}
Since $J_i(X_1,\dots,X_n:D)$ is obtained in 
general by projecting $J_i(X_1,\dots,X_n:B)$
onto $L^2(X_1,\dots,X_n,D)$, the statement $(i)$ can also be reformulated as
$$(i')\qquad J_i(X_1,\dots,X_n:D)=
J_i(X_1,\dots,X_n:B)\qquad\text{for all $i=1,\dots,n$.}$$ 
Since the implication $(ii)\Rightarrow (i)$ follows from \cite{shlyakht:cpentropy}, 
we only have 
to consider the implication $(i)\Rightarrow (ii)$.
Let us assume $(i')$. 
By \cite{V6}, we know that
\begin{align*}
j_\CC(\CC\langle X_1,\dots,X_n\rangle:D)&=\sum_{i=1}^n[J_i(X_1,\dots,X_n:D),X_i]\\
&=\sum_{i=1}^n[J_i(X_1,\dots,X_n:B),X_i]\\
&=j_\CC(\CC\langle X_1,\dots,X_n\rangle:B).
\end{align*}
Thus Theorem \ref{main2a} gives the assertion.
\end{proof}

In the same way as for \cite{V6}, Prop.5.18, we get the following consequence.

\begin{cor}
Let $(M,\tau)$ be a tracial $W^*$-probability space. Let $D\subset M$ 
be a unital $*$-subalgebra with conditional expectation $E:M\to D$ such 
that $\tau\circ E=\tau$.
Consider 
selfadjoint random variables $X_1,\dots,X_m,Y_1,\dots,Y_n\in M$ with
$\Phi^*(X_1,\dots,X_m,Y_1,\dots Y_n:D)<\infty$. 
Assume that $\dim D<\infty$.
Then the following two statements
are equivalent: 
\begin{enumerate}
\item
$\Phi^*(X_1,\dots,X_m,Y_1,\dots,Y_n:D)=\Phi^*(X_1,\dots,X_m:D)+
\Phi^*(Y_1,\dots,Y_n:D)$
\item
$\{X_1,\dots,X_m\}$ and $\{Y_1,\dots,Y_n\}$ are free with amalgamation over $D$.
\end{enumerate}
\end{cor}

\subsection{Maximization of free entropy}
Recall that, given a $D$-probability space $(M,E_D:M\to D)$,
the $D$-valued distribution of $n$ random variables $X_1,\dots,X_n\in M$
is given by the moment series of $X_1,\dots,X_n$, i.e., it consists 
of the collection of all $D$-valued moments
$E_D(d_0X_{i_1}d_1\cdots X_{i_{m-1}}d_{m-1}X_{i_m}d_m)$ for all
$m\geq 1$, all $d_1,\dots,d_m\in D$, and all $1\leq i_1,\dots,i_m\leq n$.
Since the free Fisher information $\Phi^*(X_1,\dots,X_n:D)$ with respect to
$D$ depends only on the $D$-valued distribution of $X_1,\dots,X_n$, we can 
interpret Theorem \ref{main2} also as a result about the minimization of free
Fisher informations in the following form.

\begin{cor}
\label{max:Fisher}
Let $(M,\tau)$ be a tracial $W^*$-probability space
and $B\subset M$ a unital $*$-subalgebra. Consider 
selfadjoint random variables $X_1,\dots,X_n\in M$ 
with $\Phi^*(X_1,\dots,X_n:B)<\infty$.
Let $D\subset B$ be a unital $*$-subalgebra of $B$ with conditional
expectation $E:M\to D$ such that $\tau\circ E=\tau$.
Assume that $\dim D<\infty$.
Then the following two statements are equivalent:
\begin{enumerate} 
\item
$\Phi^*(X_1,\dots,X_n:B)$ is minimal among 
$\Phi^*(Y_1,\dots,Y_n:B)$ for all $B$-valued random variables
$Y_1,\dots,Y_n$ which have the same $D$-valued distribution as
$X_1,\dots,X_n$.
\item
$\{X_1,\dots,X_n\}$ is free from $B$ with amalgamation over $D$.
\end{enumerate}
\end{cor}

It is even more striking to formulate this statement in a dual
version as a maximization result for free entropy. Let us first recall that --
for a tracial $W^*$-probability space $(M,\tau)$, a unital $*$-subalgebra 
$B\subset M$, and selfadjoint random variables $X_1,\dots,X_n\in M$ --  
the 
\emph{relative free entropy $\chi^*(X_1,\dots,X_n:B)$ with respect to $B$} is defined
as (see \cite{V5})
\begin{multline}
\chi^*(X_1,\dots,X_n:B):=\\
\frac 12\int_0^\infty\bigl(\frac n{1+t}-
\Phi^*(X_1+\sqrt tS_1,\dots,X_n+\sqrt tS_n:B)\bigr)dt +\frac n2 \log(2\pi e),
\end{multline}
where the $S_j$'s are $(0,1)$-semicircular and $B\langle X_1,\dots,X_n\rangle,
\{S_1\},\dots,\{S_n\}$ are free.
 
\begin{thm}
\label{entropy}
Let $(M,\tau)$ be a tracial $W^*$-probability space
and $B\subset M$ a unital $*$-subalgebra. Consider 
selfadjoint random variables $X_1,\dots,X_n\in M$ 
with $\chi^*(X_1,\dots,X_n:B)>-\infty$.
Let $D\subset B$ be a unital $*$-subalgebra of $B$ with conditional
expectation $E:M\to D$ such that $\tau\circ E=\tau$.
Assume that $\dim D<\infty$.
Then the following two statements are equivalent:
\begin{enumerate} 
\item
$\chi^*(X_1,\dots,X_n:B)=\chi^*(X_1,\dots,X_n:D)$.
\item
$\{X_1,\dots,X_n\}$ is free from $B$ with amalgamation over $D$.
\end{enumerate}
\end{thm}

\begin{proof}
That the freeness condition implies equality of the entropies is
contained in \cite{shlyakht:cpentropy}, so we only have to prove
the other implication. One should note that one can use the same
set of semicirculars $S_1,\dots,S_n$ for dealing with the free
entropy with respect to $B$ as well as with the free entropy with respect
to $D$. But then the assumption $(i)$ implies that we have the
equality 
$$\Phi^*(X_1+\sqrt t S_1,\dots,X_n+\sqrt t S_n:B)=
\Phi^*(X_1+\sqrt t S_1,\dots,X_n+\sqrt t S_n:D)$$ 
for almost all $t>0$.
By our corresponding result for the free Fisher 
information, Theorem \ref{main2}, this implies that
$\{X_1+\sqrt t S_1,\dots,X_n+\sqrt t S_n\}$ is free from $B$ with
amalgamation over $D$ for almost all $t>0$. By 
letting $t\rightarrow 0$ and by using the continuity of
the moments in $t$, we get the assertion. 
\end{proof}

\begin{cor}
\label{max:entropy}
Let $(M,\tau)$ be a tracial $W^*$-probability space
and $B\subset M$ a unital $*$-subalgebra. Consider 
selfadjoint random variables $X_1,\dots,X_n\in M$ 
with $\chi^*(X_1,\dots,X_n:B)>-\infty$.
Let $D\subset B$ be a unital $*$-subalgebra of $B$ with conditional
expectation $E:M\to D$ such that $\tau\circ E=\tau$.
Assume that $\dim D<\infty$.
Then the following two statements are equivalent:
\begin{enumerate} 
\item
$\chi^*(X_1,\dots,X_n:B)$ is maximal among 
$\chi^*(Y_1,\dots,Y_n:B)$ for all $B$-valued random variables
$Y_1,\dots,Y_n$ which have the same $D$-valued distribution as
$X_1,\dots,X_n$.
\item
$\{X_1,\dots,X_n\}$ is free from $B$ with amalgamation over $D$.
\end{enumerate}
\end{cor}

\subsection{$D\subset B$-Haar unitary elements}
It is to be expected (but still unproven)
that in general the definition for $\chi^*$ agrees with the one coming
from the micro-states approach and thus $\chi^*$ should be a measure of the randomness
of the given distribution. In that interpretation the
above statement is quite plausible: maximal randomness for the $B$-valued
distribution under the constraint of a fixed $D$-valued distribution is achieved
by making the variables as free as possible from $B$ modulo the given constraint
(i.e. free with amalgamation over $D$).

One should note that fixing the $D$-valued distribution and requiring
that the variables are free from $B$ with amalgamation over $D$ determines
uniquely the $B$-valued distribution. Thus for each $D$-valued distribution there 
is exactly one $B$-valued
distribution with maximal entropy. 
In this context the question arises whether there is an explicit way of
realizing the situation with maximal entropy. A universal way of doing so is
by conjugating with a $D\subset B$-version of a Haar unitary

\begin{defn}
A unitary $u$ is called a \emph{$D\subset B$-Haar unitary element}, if the following
requirements are fulfilled:
\begin{enumerate}
\item
$\{u,u^*\}$ commutes with $D$
\item
$E_D(u^k)=\delta_{k0}$ for all $k\in\ZZ$
\item
$\{u,u^*\}$ is free from $B$ with amalgamation over $D$
\end{enumerate}
\end{defn}

The last two requirements determine the $B$-valued distribution of $u,u^*$
uniquely, and the first condition has the effect that conjugating with $u$
does not change the $D$-valued distribution. 

\begin{prop}
Consider a $B$-valued
random variable $X$ and choose a $D\subset B$-Haar unitary $u$ such that
$\{u,u^*\}$ is free from $X$ with amalgamation over $B$. Put $Y:=uXu^*$.
Then the $D$-valued distribution of $Y$ is the same as the $D$-valued 
distribution of $X$, but $Y$ is free from $B$ 
with amalgamation over $D$, i.e.,
$Y$ maximizes the free entropy with respect to $B$ under all variables which
have the same $D$-valued distribution as $X$. 
\end{prop}

\begin{proof}
We have to check that $Y$ is free from $B$ with amalgamation over $D$. 
By our assumptions, we have that $X$ is free from $\{u,u^*\}$ with 
amalgamation over $B$ and 
that $\{u,u^*\}$ is free from $B$ with amalgamation over D. 
But then a slightly modified version of our
Corollary
\ref{cor:free} implies that $B\langle X\rangle$ is free from $\{u,u^*\}$
with amalgamation over D. Then it follows directly from the definition of
freeness that $uXu^*$ is free from $B$ with amalgamation over $D$.
\end{proof}

In some sense,
conjugating with $u$ can be considered as a random rotation of the degrees
of freedom of $X$ which are not fixed by the $D$-valued distribution. 
Similar constructions are possible
for more than one variable.

\subsection{Example: $R$-cyclic matrices}
The motivating example for our investigations on equality of Fisher
informations for different subalgebras was the following special case:
$M=M_d(A)=M_d(\CC)\otimes A$, $B=M_d(\CC)$, and 
$D\subset B$ is the unital $*$-algebra of constant diagonal matrices, i.e.
$$D=\{\begin{pmatrix}
\alpha_1&\dots&0\\
\vdots&\ddots&\vdots\\
0&\dots&\alpha_d
\end{pmatrix}\mid
\alpha_1,\dots,\alpha_d\in\CC\}.$$ 
In this case, our random variables $X\in M$ are $d\times d$-matrices with
entries from $A$,
$$X=(x_{ij})_{i,j=1}^d\qquad\text{with}\qquad x_{ij}\in A,$$
and statements about operator-valued properties of $X$ can also be reformulated
as scalar-valued properties of the entries $x_{ij}$.

It is quite instructive to see that in this case
the $B$-valued distribution of $X=(x_{ij})$ is
the same as the joint distribution of all entries (i.e., the collection of all
possible moments of the $x_{ij}$), whereas the $D$-valued distribution of $X$  is
given by 
the collection of all cyclic moments of the $x_{ij}$, i.e. by all
moments of the form
$\tau(x_{i_1i_2}x_{i_2i_3}\cdots x_{i_ni_1})$
for all integer $n$ and all $1\leq i_1,\dots,i_n\leq d$.

In \cite{NSS2}, we showed that the statement $X=(x_{ij})$ is
free from $B$ with amalgamation over $D$ is equivalent to the fact that the
family $\{x_{ij}\mid i,j=1\dots d\}$ is $R$-cyclic, 
which means the following: all cumulants
$\kk_n(x_{i_1j_1},x_{i_2j_2},\dots,x_{i_nj_n})$ vanish for which it is not
true that $j_1=i_2,j_2=i_3,\dots,j_n=i_i$. 

Furthermore, in \cite{NSS1}, we showed that the operator-valued free 
Fisher information $\Phi^*(X:M_d(\CC))$
of the matrix $X$ with respect to $M_d(\CC)$ is, up to a factor $d^3$, the same as the
scalar-valued free Fisher information $\Phi^*(x_{ij}\mid i,j=1,\dots,d)$ of
the entries of the matrix $X$ (where we used a slight extension of the definition
of $\Phi^*$ to the case where some of the arguments are not self-adjoint itself,
but come always in pairs with their adjoint). 

Thus, in this special case, 
we can rewrite Corollary \ref{max:entropy} from operator-valued properties of the
matrix $X=(x_{ij})$ 
to a form which
involves only scalar-valued properties of the entries $x_{ij}$.
Of course, a similar version holds for minimization of $\Phi^*$ instead of 
maximization of $\chi^*$.
 
\begin{cor}
Let $(A,\tau)$ be a tracial $W^*$-probability space and consider
random variables $x_{ij}\in A$ ($i,j=1,\dots,d$) with 
$\chi^*(\{x_{ij}\}_{i,j=1}^d)>-\infty$. Then the following statements are
equivalent:
\begin{enumerate}
\item
$\chi^*(\{x_{ij}\}_{i,j=1}^d)$ is minimal among  
$\chi^*(\{y_{ij}\}_{i,j=1}^d)$ for all $\{y_{ij}\}_{i,j=1}^d$ which have the same cyclic
moments as $\{x_{ij}\}_{i,j=1}^d$.
\item
The family $\{x_{ij}\}_{i,j=1}^d$ 
is $R$-cyclic.
\end{enumerate}
\end{cor}

One should also note that in this case a $D\subset B$-Haar unitary element $u$
is given by a diagonal matrix whose non-vanishing entries are free Haar unitaries.
The condition that such a $u$ is free from $X$ with amalgamation over $B$ just means
that all entries of $u$ are free from all entries of $X$.

\providecommand{\bysame}{\leavevmode\hbox to3em{\hrulefill}\thinspace}

\end{document}